\documentclass[a4,11pt]{article}
\usepackage{amsmath,amssymb,amsthm}
\usepackage[arrow,matrix]{xy}
\fontsize{20}{8}\selectfont
\theoremstyle{remark}

\theoremstyle{definition}

\theoremstyle{question}

\theoremstyle{corollary}

\theoremstyle{example}

\newcommand{\al}{\alpha}

\begin{document}

\begin{center}
{\Large {\bf CRITERIA FOR THE BOUNDEDNESS OF POTENTIAL OPERATORS IN GRAND LEBESGUE SPACES}}
\end{center}
\vskip+1cm
\begin{center}
{\bf ALEXANDER  MESKHI}
\end{center}




\vskip+2cm

{\bf Abstract.} It is shown that that the fractional integral
operators with the parameter $\alpha$, $0<\alpha<1$,  are not bounded between the
generalized grand Lebesgue spaces $L^{p), \theta_1}$ and  $L^{q),
\theta_2}$ for $\theta_2 < (1+\alpha q)\theta_1$, where
$1<p<1/\alpha$ and $q=\frac{p}{1-\alpha p}$.  Besides this, it is proved that the one--weight inequality
$$
\|I_{\alpha}(fw^{\alpha})\|_{L_{w}^{q),\theta(1+\alpha
q)}}\leq c\|f\|_{L_{w}^{p),\theta}},
$$
where $I_{\alpha}$ is the Riesz potential operator on the interval
$[0,1]$, holds if and only if $w\in A_{1+q/p'}$.

\vskip+1cm

{\bf 2010 Mathematics Subject Classification}: 42B25, 46E30.
\vskip+0.4cm

{\bf Key words and phrases:} Grand Lebesgue spaces, potential operators, weights, Boundedness, one--weight inequality.

\vskip+2cm

\section*{Introduction}

In this paper we show that potential operators with the parameter
$\alpha$, $0<\alpha<1$, are not bounded from $L^{p)}$ to $L^{q)}$,
where $1<p<\infty$ and $q$ is the Hardy--Littlewood--Sobolev
exponent of $p$: $q= \frac{p}{1-\alpha p}$.   This phenomena
motivates us to investigate the boundedness problem for the Riesz
potential operator $I_{\alpha}$   in the generalized  grand
Lebesgue spaces. In particular, we study this problem in  $L^{p,
\theta}_w$ spaces and prove that the one--weight inequality
$$
\|I_{\alpha}(fw^{\alpha})\|_{L_{w}^{q),\theta(1+\alpha
q)}([0,1])}\leq c\|f\|_{L_{w}^{p),\theta}([0,1])}
$$
holds if and only if $w$ belongs to the Muckenhoupt's class
$A_{1+q/p'}$.

The unweight spaces $L^{p, \theta}$ (i.e. $L^{p, \theta}_w$ for $w\equiv
const$) were introduced by E. Greco, T. Iwaniec and C. Sbordone
\cite{GrIwSb} when they studied existence and uniqueness of the nonhomogeneous $n-$
harmonic equation $div A(x, \nabla u)= \mu$.

The grand Lebesgue spaces $L^{p)}= L^{p), 1}$ first appeared in
the paper by T. Iwaniec and C. Sbordone \cite{IwSb}. In that paper
the authors showed that if $f= (f_1, \cdots, f_n): \Omega \to
{\Bbb{R}}^n$ belongs to the Sobolev class $W^{1,1}$, where $\Omega$ is an open subset in ${\Bbb{R}}^n$,
$n\geq 2$, then the Jacobian determinant  $J= J(f,x)= det \;Df(x)$
($J(x,f)\geq 0$ a.e.) of $f$ belongs to the class
$L^1_{loc}(\Omega)$ provided that $ g \in L^{n)}$, where
$$ g(x):= |Df(x)|= \{ \sup |Df(x) y|: y \in S^{n-1} \}. $$

Recently necessary and sufficient conditions guaranteeing  the one--weight
inequality for the Hardy--Littlewood maximal operator in
$L^{p)}_w(I)$, where $I=[0,1]$, were established by A.
Fiorenza, B. Gupra and P. Jain \cite{FiGuJa}, while the same
problem for the Hilbert transform was studied in the paper
\cite{KoMe}. In particular, it turned out that the
Hardy--Littlewood maximal operator (resp. the Hilbert transform)
is bounded in $L^{p)}_w(I)$ if and only if the weight $w$ belongs
to the Muckenhoupt class $A_p(I)$.

\section{Preliminaries}

Let $\Omega$ be a bounded subset of ${\Bbb{R}}^n$ and let $w$ be
an a.e. positive, integrable function on $\Omega$ (i.e. a weight).
The weighted generalized grand Lebesgue space $L^{p),
\theta}(\Omega)$ ($1<p<\infty$) is the class of those $f:
\Omega\to {\Bbb{R}}$ for which  the norm
$$  \|f\|_{L^{p), \theta}_w(\Omega)} = \sup_{0<\varepsilon <p-1} \bigg(\frac{\varepsilon^{\theta}}{|\Omega|} \int_{\Omega}
|f(t)|^{p-\varepsilon} w(t) dt \bigg)^{1/(p-\varepsilon)} $$ is
finite.


If $w\equiv 1$, then we denote $L^{p), \theta}(\Omega):=L^{p),
\theta}_1(\Omega)$. The space $L^{p), \theta}_w(\Omega)$ is not
rearrangement invariant unless $w\equiv \; const$.

H\"older's inequality and simple
estimates yield the following embeddings (see also \cite{GrIwSb},
\cite{FiGuJa}):

$$ L_{w}^{p}(\Omega)\subset L_{w}^{p),\theta_{1}}(\Omega)\subset  L_{w}^{p),\theta_{2}}(\Omega)\subset
L_{w}^{p-\varepsilon}(\Omega),\eqno{(1.1)} $$ where $0<
\varepsilon< p-1$ and $\theta_1 < \theta_2. $

In the classical weighted Lebesgue spaces $L^p_w$ the equality

$$ \|f\|_{L^{p}_w}= \|w^{1/p}f \|_{L^{p}}$$
holds but this property fails in the case of grand Lebesgue
spaces. In particular, there is $f\in L^{p)}_w$ such that $w^{1/p}
f \notin L^{p)}$ (see also \cite{FiGuJa} for the details).

Let $\varphi$ be positive increasing function on $(0,p-1)$
satisfying the condition $\varphi(0+)=0$, where $1<p<\infty$. We
will also need the following auxiliary class of functions defined on
$\Omega$ and associated with $\varphi$:
$$
L_{w}^{p),\varphi(x)}(\Omega):=\left\{f:\;\sup\limits_{0<\varepsilon\leq
p-1}\left(\varphi(\varepsilon)^{\frac{1}{p-\varepsilon}}\|f\|_{L_{w}^{p-\varepsilon}}\right)<\infty\right\}.
$$

The space $L^{p), \theta}_w(\Omega)$, $\theta>0$, is the special case of $L_{w}^{p),\varphi(x)}(\Omega)$ taking $\varphi (x) = \frac{x^{\theta}}{|\Omega|}$.


Throughout the paper the symbol $\varphi(t)\approx\psi(t)$ means
that there exist positive constants $c_{1}$ and $c_{2}$ such that
$c_{1}\varphi(t)\leq\psi(t)\leq c_{2}\psi(t)$. Constants (often
different constants in the same series of inequalities) will
generally be denoted by $c$ or $C$. By the symbol $p'$ we denote
the conjugate number of $p$, i.e. $p' := \frac{p}{p-1}$,
$1<p<\infty$.

\section{Fractional Integrals and Fractional Maximal Functions in Unweighted Grand Lebesgue Spaces}

Let
$$
(I_\al f)(x)=\int_0^1 \frac{f(y)}{|x-y|^{1-\al}}\,dy,\;\;\;0<\al<1
$$
be the Riesz potential operator defined on $[0,1]$.  We begin with
the following result: \vskip+0.2cm

{\bf Theorem 2.1.} {\em  Let $0<\alpha<1$, $1<p<\frac{1}{\alpha}$,
$\theta_{1}$ and $\theta_{2}$ be positive numbers such that
$\theta_{2}<\theta_{1}(1+\alpha q)$, where $q=\frac{p}{1-\alpha
p}$. Then the operator $I_{\alpha}$ is not bounded from
$L^{p),\theta_{1}}$ to $L^{q),\theta_{2}}$.}

\begin{proof}
Suppose the contrary: $I_{\alpha}$ is bounded from
$L^{p),\theta_{1}}$ to $L^{q),\theta_{2}}$ i. e. the inequality
$$
\|I_{\alpha}f\|_{L^{q),\theta_{2}}}\leq c\|f\|_{L^{p),\theta_{1}}} \eqno{(2.1)}
$$
holds, where the positive constant $c$ does not depend on $f$.
Taking $f=\chi_{J}$ in (2.1), where $J$ is an interval in $[0,1]$,
we have
$$
(I_{\alpha}f)(x)=\int\limits_{J}\frac{dy}{|x-y|^{1-\alpha}}\geq
|J|^{\alpha},\;\;x\in J.
$$
Consequently,
$$
\|I_{\alpha}f\|_{L^{q),\theta_{2}}}\geq
|J|^{\alpha}\|\chi_{J}\|_{L^{q),\theta_{2}}}.
$$
Taking  inequality (2.1) into account we have that
$$
|J|^{\alpha}\|\chi_{J}\|_{L^{q),\theta_{2}}}\leq
c\|\chi_{J}\|_{L^{p),\theta_{1}}}, \eqno{(2.2)}
$$
where the positive constant $c$ does not depend on $J$.

Let us define the number $\varepsilon_{J}$ which is between $0$
and $p-1$ and satisfies the condition
$$
\sup\limits_{0<\varepsilon\leq
p-1}\left(\varepsilon^{\theta_{1}}|J|\right)^{\frac{1}{p-\varepsilon}}=
\left(\varepsilon_{J}^{\theta_{1}}|J|\right)^{\frac{1}{p-\varepsilon_{J}}}
\eqno{(2.3)}
$$
Now we claim that $\lim\limits_{|J|\rightarrow
0}\varepsilon_{J}=0$. Indeed, suppose the contrary: that there is a sequence of
intervals $J_{n}$ and a positive number $\lambda$ such that $|J_{n}|\rightarrow 0$ and
$\varepsilon_{J_{n}}\geq\lambda>0$ for all $n\in N$. It is obvious that  we
can choose $J_{n_{0}}$ so that
$$
\frac{|J_{n_{0}}|^{\frac{1}{\theta_{1}}}(p-1)}{e}<e^{-\frac{p}{\lambda/2}}.$$
Now we claim that $f'(x)<0$ for all $x\in [\lambda/2, p-1]$, where
$f(x)=\left(x^{\theta_{1}}|J_{n_{0}}|\right)^{\frac{1}{p-x}}$.
Indeed, it is easy to see that for $\lambda/2\leq x \leq p-1$, the
inequalities

$$\frac{|J_{n_{0}}|^{ \frac{1}{\theta_{1}}} x}{e} \leq
\frac{|J_{n_{0}}|^{\frac{1}{\theta_{1}}}(p-1)}{e}<e^{-\frac{p}{\lambda/2}}
\leq e^{-\frac{p}{x}}.  $$
hold. Hence, using the formula
$$
f'(x)=f(x)\cdot\frac{1}{p-x}\left[\frac{ln\left(x^{\theta_{1}}|J_{n_{0}}|\right)}{p-x}+\frac{\theta_{1}}{x}\right]
$$
and the fact that
$$
f'(x)<0\Longleftrightarrow
\frac{x|J_{n_{0}}|^{\frac{1}{\theta_{1}}}}{e}<e^{-\frac{p}{x}}
$$ we conclude that $f'(x)<0$.

This observation together with the equality
$\lim\limits_{x\rightarrow 0}f(x)=0$ gives that
$\varepsilon_{J_{n_{0}}}<\lambda$, where $\varepsilon_{J_{n_{0}}}$
is defined by
$$
\sup\limits_{0<\varepsilon\leq
p-1}\left(\varepsilon^{\theta_{1}}|J_{n_{0}}|\right)^{\frac{1}{p-\varepsilon}}=
\left(\varepsilon_{J_{n_{0}}}^{\theta_{1}}|J_{n_{0}}|\right)^{1/(p-\varepsilon_{J_{n_{0}}})}.
$$
This contradicts the assumption that
$\varepsilon_{J_{n}}\geq\lambda>0$ for all $n$. Further, we choose
$\eta_{J}$ so that
$$
\alpha=\frac{1}{p}-\frac{1}{q}=\frac{1}{p-\varepsilon_{J}}-\frac{1}{q-\eta_{J}}.
$$
This is equivalent to say that
$$
\eta_{J}=q-\frac{p-\varepsilon_{J}}{1-\alpha(p-\varepsilon_{J})}.
\eqno{(2.4)}
$$
By (2.2) and (2.3) we have that
$$
|J|^{\alpha}\eta_{J}^{\frac{\theta_{2}}{q-\eta_{J}}}|J|^{\frac{1}{q-\eta_{J}}}\leq
c\varepsilon_{J}^{\frac{\theta_{1}}{p-\varepsilon_{J}}}|J|^{\frac{1}{p-\varepsilon_{J}}}.
\eqno{(2.5)}
$$
(here we used the fact that if $\varepsilon_{J}$ is small, then
$0<\eta_{J}<q-1$). Now (2.5) yield:
$$
\eta_{J}^{\frac{\theta_{2}}{q-\eta_{J}}}\varepsilon_{J}^{-\frac{\theta_{1}}{p-\varepsilon_{J}}}\leq
c. \eqno{(2.6)}
$$
Further, (2.4) and (2.6) imply
$$
\left(\frac{q-\frac{p-\varepsilon_{J}}{1-\alpha(p-\varepsilon_{J})}}
{\varepsilon_{J}}\right)^{\frac{\theta_{2}}{p-\varepsilon_{J}}-\alpha\theta_{2}}
\varepsilon_{J}^{-\frac{\theta_{1}}{p-\varepsilon_{J}}+\frac{\theta_{2}}{p-\varepsilon_{J}}-\alpha\theta_{2}}
\leq c. \eqno{(2.7)}
$$
Passing now to the limit as $|J|\rightarrow 0$ we see that the
left-hand side of (2.7) tends to $+\infty$ because the limit of
the first factor is $\left[\frac{1}{(1-\alpha
p)^{2}}\right]^{\frac{\theta_{2}}{p}-\alpha\theta_{2}},$ and
$$
\lim\limits_{|J|\rightarrow
0}\varepsilon_{J}^{\frac{\theta_{2}-\theta_{1}}{p-\varepsilon_{J}}-\alpha\theta_{2}}
=\lim\limits_{|J|\rightarrow
0}\varepsilon_{J}^{\frac{\theta_{2}-\theta_{1}}{p}-\alpha\theta_{2}}=\infty
$$
(Here we used the observation
$\frac{\theta_{2}}{\theta_{1}}<1+\alpha
q\Longleftrightarrow\frac{\theta_{2}-\theta_{1}}{p}-\alpha\theta_{2}<0$).
\end{proof}

\vskip+0.2cm

Analysing the proof of Theorem 2.1 we have the result similar to that of the previous statement for the fractional maximal operator
$$M_{\alpha}f(x)=\sup_{\substack{J\ni x \\
J\subset[0,1]}}\frac{1}{|J|^{1-\alpha}}\int\limits_{J}|f|, \;\;\;
x\in [0,1].$$

\vskip+0.2cm

{\bf Theorem 2.2.} {\em  Let the conditions of Theorem $2.1$ be satisfied.  Then the operator
$M_{\alpha}$ is not
bounded from $L^{p),\theta_{1}}$ to $L^{q),\theta_{2}}$.}

\begin{proof}
Proof is the same as in the case of Theorem 2.1. We only need to
observe that the inequality
$$
M_{\alpha}f(x)\geq\frac{1}{|J|^{1-\alpha}}\int\limits_{J}dx=|J|^{\alpha},\;\;x\in
J,
$$
holds for $f(x)=\chi_{J}(x)$, where $J$ is a subinterval of $[0,1]$. Details are omitted.
\end{proof}


\section{Sobolev's Embedding in Weighted Generalized Grand Lebesgue
Spaces}

\
This section is devoted to the investigation of the one--weight inequality for the operator $I_{\alpha}$ in $L^{p), \theta}_w$ spaces.

First we introduce the function
$$
\varphi(u)=\left[\frac{u-q}{1-\alpha(u-q)}+p\right]^{1-
(u-q)\alpha}\eqno{(3.1)}
$$
where $0<\alpha<1$, $1<p<\frac{1}{\alpha}$, $q=\frac{p}{1-\alpha
p}$.
\vskip+0.2cm

To prove the main results we need some auxiliary statements.
\vskip+0.2cm

{\bf Lemma 3.1.} {\em  $\varphi(x)\approx x^{1+\alpha q}$ near
$0$.} \vskip+0.2cm The proof is straightforward and therefore is
omitted.

\vskip+0.2cm

{\bf Lemma 3.2.} {\em  Let $1<q<\infty$ and let $w$ be a weight.
Then
$$
\|f\|_{L_{w}^{q),\varphi(x)}([0,1])}\approx\|f\|_{L_{w}^{q),1+\alpha q}([0,1])}
$$ where $\varphi$ is defined by $(3.1)$.}

\vskip+0.2cm

\begin{proof}
Follows immediately from Lemma 3.1.
\end{proof}
\vskip+0.2cm
 {\bf Lemma 3.3.}  {\em Let $1<q<\infty$ and let $\theta>0$. Then
$$ \|f\|_{L_{w}^{q),\varphi(x^{\theta})}([0,1])}\approx\|f\|_{L_{w}^{q),\theta(1+\alpha
q)}([0,1])}, $$
where $\varphi$ is defined by $(3.1)$}

\vskip+0.2cm

The proof follows immediately from Lemma 3.1.

\vskip+0.2cm

{\bf Lemma 3.4.} {\em   Let $1<p<\infty$ and let $\varphi$ be as
above. Then there is a positive constant $c$ such that for all
intervals $J\subset [0,1]$ and $f\in L_{w}^{p,\varphi(x)}$ the
inequality
$$
\|f\|_{L_{w}^{p),\varphi(x)}(J)}\leq
c(w(J))^{-\frac{1}{p}}\left(\int\limits_{J}|f(t)|^{p}w(t)dt\right)^{\frac{1}{p}}\|\chi_{J}\|_{L_{w}^{p),\varphi(x)}}
$$
holds. } \vskip+0.2cm

\begin{proof}
We have
$$
\|f\|_{L_{w}^{p),\varphi(x)}(J)}=\sup\limits_{0<\varepsilon\leq
p-1}\left(\varphi(\varepsilon)\int\limits_{J}
|f(x)|^{p-\varepsilon}w(x)dx\right)^{\frac{1}{p-\varepsilon}}
$$

$$
=\sup\limits_{0<\varepsilon\leq
p-1}\left(\varphi(\varepsilon)\int\limits_{J}
|f(x)|^{p-\varepsilon}w(x)^{\frac{p-\varepsilon}{p}}w(x)^{\frac{\varepsilon}{p}}dx\right)^{\frac{1}{p-\varepsilon}}
$$

$$
\leq\sup\limits_{0<\varepsilon\leq
p-1}\varphi(\varepsilon)^{\frac{1}{p-\varepsilon}}\left(\int\limits_{J}
\left(|f(x)|^{p-\varepsilon}w(x)^{\frac{p-\varepsilon}{p}}\right)^{\frac{p}{p-\varepsilon}}dx\right)^{\frac{1}{p}}
\left(\int\limits_{J}\left[w^{\frac{\varepsilon}{p}}(x)\right]^{\frac{p}{\varepsilon}}dx\right)^
{\frac{\varepsilon}{p(p-\varepsilon)}}
$$

$$
=\sup\limits_{0<\varepsilon\leq
p-1}\varphi(\varepsilon)^{\frac{1}{p-\varepsilon}}\left(\int\limits_{J}
|f(x)|^{p}w(x)dx\right)^{\frac{1}{p}}\left(\int\limits_{J}w(x)dx\right)^
{\frac{\varepsilon}{p(p-\varepsilon)}}
$$

$$
=\left(\int\limits_{J}
|f(x)|^{p}w(x)dx\right)^{\frac{1}{p}}\left(\int\limits_{J}w(x)dx\right)^{-\frac{1}{p}}
\sup\limits_{0<\varepsilon\leq
p-1}\left(\varphi(\varepsilon)\int\limits_{J}w(x)dx\right)^{\frac{1}{p-\varepsilon}}
$$

$$
=\left(\int\limits_{J}|f(x)|^{p}w(x)dx\right)^{\frac{1}{p}}(w(J))^{-\frac{1}{p}}\|\chi_{J}\|_{L_{w}^{p),\varphi(x)}(J)}.
$$
\end{proof}

{\bf Lemma 3.5.}  {\em Let $\theta >0$, $1<p<\infty$, $0<\alpha<1/p$ and let
$q=\frac{p}{1-\alpha p}$. Suppose that the inequality
$$
\|I_{\alpha}(fw^{\alpha})\|_{L_{w}^{q),\theta}([0,1])}\leq
c\|f\|_{L_{w}^{p),\theta}([0,1])} \eqno{(3.2)}
$$
holds. Then
$$
\int\limits_{0}^{1}w^{-p'/q}(x)dx<\infty.
$$}

\begin{proof}
Suppose the contrary:
$\int\limits_{0}^{1}w^{-p'/q}(x)dx=\|w^{\alpha-1}\|_{L_{w}^{p'}}=\infty$.
This means that there is a function $g\in L_{w}^{p}$ such that
$\int\limits_{0}^{1}g w^{\alpha}=\infty$.

On the other hand,
$$
I_{\alpha}(g w^{\alpha})(x)=\int\limits_{0}^{1}\frac{g(t)w^{\alpha}(t)}{|x-t|^{1-\alpha}}dt\geq\int\limits_{0}^{1}
g(t)w^{\alpha}(t)dt=\infty, \;\;\; x\in [0,1].
$$
Further, Lemma 3.4 with $\varphi(x) = x^{\theta}$ implies that $g\in L_{w}^{p),\theta}([0,1])$. But
$I_{\alpha}(gw^{\alpha})(x)=\infty$ for $x\in [0,1]$. This contradicts inequality
(3.2).
\end{proof}

\vskip+0.1cm

{\bf Definition 3.1.} Let $1<r<\infty$.  We say that a weight function $w$ belongs to the Muckenhoupt's class $A_{r}([0,1])$ ($w\in A_{r}([0,1])$) if
$$ A_{r}(w) :=  \sup\limits_{J\subset[0,1]}\left(\frac{1}{|J|}\int\limits_{J}w\right)^{1/r}\left(\frac{1}{|J|}\int\limits_{J}w
^{1-r'}\right)^{1/r'}<\infty,
$$
where the supremum is taken over all subintervals $J$ of $[0,1]$.
\vskip+0.2cm

{\bf Lemma 3.6.} {\em  Let $0<\alpha<1$, $1<p<1/\alpha$. We set $q= \frac{p}{1-\alpha p}$. Suppose that $w\in A_{1+q/p'}([0,1])$, i.e.,

$$ \sup\limits_{J\subset[0,1]}\left(\frac{1}{|J|}\int\limits_{J}w\right)^{1/q}\left(\frac{1}{|J|}\int\limits_{J}w
^{-p'/q}\right)^{1/p'}<\infty. $$

Then there are positive constants $\sigma_1$, $\sigma_2$ and $L$ satisfying the conditions:

$$\frac{1}{p-\sigma_2} - \frac{1}{q-\sigma_1}=\alpha, \;\; w\in A_{1+\frac{q- \sigma_1}{(p-\sigma_2)'}},$$
$$ \| K_{\alpha} \|_{L^{p-\eta}_w \to L^{q-\varepsilon}_w} \leq L $$
for all $0\leq  \varepsilon\leq \sigma_1$, $0\leq \eta\leq
\sigma_2$ with $ \frac{1}{p-\eta}- \frac{1}{q-\varepsilon}
=\alpha$, where $K_{\alpha}$ is the operator defined as follows
$K_{\alpha}f= I_{\alpha}(fw^{\alpha})$.}

\vskip+0.2cm

{\em Proof.} Since $w\in A_{1+q/p'}$ by the openness property of Muckenhoupt's classes (see \cite{Mu}) we have that there are small  positive numbers $\sigma_1 $ and $\sigma_2$ such that $\frac{1}{p-\sigma_2} - \frac{1}{q-\sigma_1} = \alpha$ and $w\in A_{1+(q-\sigma_1)/(p-\sigma_2)'}$.

By the result of B. Muckenhoupt and R. L. Wheeden \cite{MuWh} we
have that the operator $K_{\alpha}$ is bounded from $L^{p}_w$ to
$L^q_w$ and from $L^{p-\sigma_2}_w$ to $L^{q-\sigma_1}_w$. Let
$0<t<1$ and let us define positive numbers $\eta$ and
$\varepsilon$ so that
$$ \frac{1}{p-\eta}= \frac{t}{p}+ \frac{1-t}{p-\sigma_2}, \;\;\;
\frac{1}{q-\varepsilon}= \frac{t}{q}+ \frac{1-t}{q-\sigma_1}. $$
Then by applying the Rieasz--Thorin theorem (see e.g. \cite{Dua},
p. 16) we have that $K_{\alpha}$ is bounded from $ L^{p-\eta}$ to
$L^{q-\varepsilon}$ and moreover,
$$\|K_{\alpha}\|_{L^{p-\eta}_w \to L^{q-\varepsilon}_w} \leq \|K_{\alpha}\|^{t}_{L^{p}_w\to L^{q}_w}
\|K_{\alpha}\|^{1-t}_{L^{p-\sigma_2}_w\to L^{q-\sigma_1}_w}. $$
Observe now that
$$ \frac{1}{p-\eta} - \frac{1}{q-\varepsilon}= \frac{t}{p}- \frac{t}{q} + \frac{1-t}{p-\sigma_2}-
\frac{1-t}{q- \sigma_1}$$

$$ =t\big( \frac{1}{p}- \frac{1}{q} \big) + (1-t) \Big( \frac{1}{p-\sigma_2}-\frac{1}{q- \sigma_1}\Big) =
t \alpha+ (1-t) \alpha = \alpha.$$

The lemma is proved since we can take $L=
\|K_{\alpha}\|_{L^{p}_w\to L^{q}_w}
\|K_{\alpha}\|_{L^{p-\sigma_2}_w\to L^{q-\sigma_1}_w}$ (since
without loss of generality we can assume that each term is greater
or equal to $1$).$\;\;\;\; \Box$

\vskip+0.2cm

{\bf Theorem 3.1.} {\em Let $1<p<\infty$ and let $0<\alpha<1/p$.
Suppose that $\theta>0$. We set $q=\frac{p}{1-\alpha p}$. Then the
inequality
$$
\|I_{\alpha}(fw^{\alpha})\|_{L_{w}^{q),\theta(1+\alpha
q)}([0,1])}\leq c\|f\|_{L_{w}^{p),\theta}([0,1])} \eqno{(3.3)}
$$
holds if and only if $w\in A_{1+q/p'}([0,1])$. }

\begin{proof}
By Lemma 3.1 we have that (3.3) is equivalent to the inequality
$$ \|I_{\alpha}(fw^{\alpha})\|_{L_{w}^{q),\psi(x)}([0,1])}\leq c\|f\|_{L_{w}^{p),\theta}([0,1])}, \eqno{(3.4)} $$
where
$$
\psi(x)=\varphi(x^{\theta}),\;\;\;\varphi(x)=\left[\frac{x-q}{1-\alpha(x-q)}+p\right]^{1-(x-1)\alpha}.
\eqno{(3.5)}
$$
{\em Necessity.} Let (3.3) and hence (3.4) hold. By Lemma 3.5 we
have that $\int\limits_{0}^{1}w^{-p'/q}<\infty$. Let us take
$f=\chi_{J}w^{-\alpha-p'/q}$. Then for $x\in J$, we get that
$$
I_{\alpha}(w^{\alpha}f)(x)\geq\frac{1}{|J|^{1-\alpha}}\int\limits_{J}w^{\alpha}f=\frac{1}{|J|^{1-\alpha}}
\int\limits_{J}w^{-p'/q}.
$$
Hence,
$$
\|I_{\alpha}(w^{\alpha}f)\|_{{L_{w}^{q),\psi(x)}([0,1])}}\geq
|J|^{\alpha-1}\left(\int\limits_{J}w^{-p'/q}\right)\|\chi_{J}\|
_{L_{w}^{q),\psi(x)}([0,1])}.
$$
Further, by Lemma 3.4 we find that
$$
|J|^{\alpha-1}\left(\int\limits_{J}w^{-p'/q}\right)\|\chi_{J}\|
_{L_{w}^{q),\psi(x)}([0,1])}
$$
$$
\leq c \|f\|_{L^{p),\theta}([0,1])}\leq c(w(J))^{-\frac{1}{p}}
\left(\int\limits_{J}|f(t)|^{p}w(t)dt\right)^{\frac{1}{p}}
\|\chi_{J}\| _{L_{w}^{p),\theta}([0,1])}
$$
$$
=c w(J)^{-\frac{1}{p}} \left(\int\limits_{J}w^{-p'/q}\right)^{1/p}
\|\chi_{J}\| _{L_{w}^{p),\theta}([0,1])}.
$$
Further, it is easy to see that there is a number $\eta_{J}$
depending on $J$ such that $0<\eta_{J}\leq p-1$ and
$$
|J|^{\alpha-1}w(J)^{\frac{1}{p}}\left(\int\limits_{J}w^{-p'/q}\right)^{\frac{1}{p'}}
\|\chi_{J}\| _{L_{w}^{q),\psi(x)}([0,1])} \leq
c\left(\eta_{J}w(J)\right)^{\frac{1}{p-\eta_{J}}}.
$$
For such  $\eta_{J}$ we choose $\varepsilon_{J}$ so that
$$
\frac{1}{p-\eta_{J}}-\frac{1}{q-\varepsilon_{J}}=\alpha.
$$
Then $0<\varepsilon_{J}\leq q-1$ and
$$
|J|^{\alpha-1}w(J)^{\frac{1}{p}-\frac{1}{p-\eta_{J}}}\eta_{J}^{-\frac{\theta}{p-\eta_{J}}}
\psi(\varepsilon_{J})^{\frac{1}{q-\varepsilon_{J}}}w(J)^{\frac{1}{q-\varepsilon_{J}}}
\left(\int\limits_{J}w^{-p'/q}\right)^{\frac{1}{p'}}\leq c. $$


Observe that by Lemma 3.1 we have that
$$ \eta_{J}^{-\frac{\theta}{p-\eta_{J}}}\psi(\varepsilon_{J})^{\frac{1}{q-\varepsilon_{J}}}
=\eta_{J}^{-\frac{\theta}{p-\eta_{J}}}\varphi\left(\varepsilon_{J}^{\theta}\right)^{\frac{1}{q-\varepsilon_{J}}}
\approx\eta_{J}^{-\frac{\theta}{p-\eta_{J}}}\varepsilon_{J}^{\frac{\theta(1+\alpha
q)}{q-\varepsilon_{J}}}=\left(\eta_{J}^{-\frac{1}{p-\eta_{J}}}\varepsilon_{J}^{\frac{1+\alpha
q}{q-\varepsilon_{J}}}\right)^{\theta}
$$
$$
\approx\left(\eta_{J}^{-\frac{1}{p-\eta_{J}}}\varphi(\varepsilon_{J})
^{\frac{1}{q-\varepsilon_{J}}}\right)^{\theta}=1
$$
and also,
$$
\frac{1}{p}-\frac{1}{p-\eta_{J}}+\frac{1}{q-\varepsilon_{J}}=\frac{1}{p}-\alpha=\frac{1}{q}.
$$
Finally, we have that
$$
|J|^{\alpha-1}w(J)^{\frac{1}{q}}\left(\int\limits_{J}w^{-p'/q}\right)^{1/p'}\leq
c.
$$
Necessity is proved.

{\em Sufficiency.} Using Lemma 3.6 we have that there are positive
constants $\sigma_1$, $\sigma_2$ and $L$ satisfying the
conditions: $\frac{1}{p-\sigma_2} - \frac{1}{q-\sigma_1}=\alpha,
\;\; w\in A_{1+\frac{q- \sigma_1}{(p-\sigma_2)'}}$, $ \|
K_{\alpha} \|_{L^{p-\eta}_w \to L^{q-\varepsilon}_w} \leq L $ for
all $0\leq \varepsilon\leq \sigma_1$, $0\leq \eta\leq \sigma_2$
with $ \frac{1}{p-\eta}- \frac{1}{q-\varepsilon} =\alpha$, where
$K_{\alpha}$ is the operator defined by $K_{\alpha}f=
I_{\alpha}(fw^{\alpha})$.

Let $\sigma$ be  a small positive number such that $\sigma< \sigma_1<q-1$ and let us fix $\varepsilon\in (\sigma, q-1]$. Then
$\frac{q-\sigma}{q-\varepsilon}>1$. By H\"{o}lder's inequality we have that
$$
\|I_{\alpha}(fw^{\alpha})\|_{L_{w}^{q-\varepsilon}([0,1])}\leq
\left(\int\limits_{0}^{1}|I_{\alpha}(fw^{\alpha})(x)|^{q-\sigma}w(x)dx\right)^{\frac{1}{q-\sigma}}
w([0,1])^{\frac{\varepsilon-\sigma}{(q-\sigma)(q-\varepsilon)}}
$$
because
$\left(\frac{q-\sigma}{q-\varepsilon}\right)'=\frac{q-\sigma}{\varepsilon-\sigma }$.

Further, the conditions $\sigma<q-1$ and $\sigma<\varepsilon<q-1$
yield
$$
0<\frac{\varepsilon-\sigma}{(q-\sigma)(q-\varepsilon)}<\frac{q-1-\sigma}{q-\sigma},\;\;\;
(q-1)\sigma^{-\frac{1}{q-\sigma}}>1.
$$
Consequently, using the well--known result by B. Muckenhoupt and
R. L. Wheeden \cite{MuWh} for the classical weighted Lebesgue
spaces:
$$
\|I_{\alpha}(fw^{\alpha})\|_{L_{w}^{q}([0,1])}\leq
c\|f\|_{L_{w}^{p}([0,1])} \Longleftrightarrow w\in A_{1+ q/p'}([0,1]), \;\; q=\frac{p}{1-\alpha p},
$$
we find that
$$
\|I_{\alpha}f\|_{L_{w}^{q),\psi(x)}([0,1])}=\max\bigg\{\sup\limits_{0<\varepsilon\leq\sigma}
\psi(\varepsilon)^{\frac{1}{q-\varepsilon}}\|I_{\alpha}(fw^{\alpha})\|_{L_{w}^{q-\varepsilon}([0,1])},$$

$$
\sup\limits_{\sigma<\varepsilon\leq q-1}\psi(\varepsilon)^{\frac{1}{q-\varepsilon}}\|I_{\alpha}(fw^{\alpha})\|
_{L^{q-\varepsilon}_w([0,1])}\bigg\}
$$

$$
\leq\max\bigg\{\sup\limits_{0<\varepsilon\leq\sigma}
\psi(\varepsilon)^{\frac{1}{q-\varepsilon}}\|I_{\alpha}(fw^{\alpha})\|_{L_{w}^{q-\varepsilon}([0,1])},$$

$$
\sup\limits_{\sigma<\varepsilon\leq q-1}\psi(\varepsilon)^{\frac{1}{q-\varepsilon}}\|I_{\alpha}(fw^{\alpha})\|
_{L_{w}^{q-\varepsilon}}w([0,1])^{\frac{\varepsilon-\sigma}{(q-\sigma)(q-\varepsilon)}}\bigg\}
$$

$$
\leq\max\left\{1,\;\sup\limits_{\sigma <\varepsilon\leq q-1}
\psi(\varepsilon)^{\frac{1}{q-\varepsilon}}\psi(\sigma)^{-\frac{1}{q-\sigma}}
w([0,1])^{\frac{\varepsilon-\sigma}{(q-\sigma)(q-\varepsilon)}}\right\}
\sup\limits_{0<\varepsilon\leq\sigma}
\psi(\varepsilon)^{\frac{1}{q-\varepsilon}}\|I_{\alpha}(fw^{\alpha})\|_{L_{w}^{q-\varepsilon}([0,1])}
$$

$$
\leq c\max\left\{1,\;\left[\sup\limits_{\sigma <\varepsilon\leq
q-1}(\psi(\varepsilon))^{\frac{1}{q-\varepsilon}}\right]\varphi(\sigma)^{-\frac{1}{q-\sigma}}
(1+w([0,1])^{\frac{q-1-\sigma}{q-\sigma}}\right\}
\sup\limits_{0<\eta \leq
\sigma_0}\eta^{\frac{\theta}{p-\eta}}\|f\|_{L_{w}^{p-\eta}([0,1])}
$$

$$
\leq c \bigg( \sup\limits_{\sigma<\varepsilon\leq
q-1}\psi(\varepsilon)^{\frac{1}{q-\varepsilon}}\bigg)
\varphi(\sigma)^{-\frac{1}{q-\sigma}}
(1+w([0,1]))^{\frac{q-1-\sigma}{q-\sigma}}\|f\|_{L_{w}^{p),\theta}([0,1])}.
$$
Here $\sigma_{0}$ is small positive number such that when
$0<\varepsilon\leq\sigma$, then $0<\eta\leq\sigma_{0}< \sigma_1< p-1$. Also, we used the estimates:
$$
\psi(\varepsilon)^{\frac{1}{q-\varepsilon}}\approx\varepsilon^{\frac{\theta(1+\alpha
q)}{q-\varepsilon}}\approx\varphi(\varepsilon)^{\frac{\theta}{q-\varepsilon}}=\eta^{\frac{\theta}{p-\eta}},\;\; \text{as}\;\;\;
\varepsilon\rightarrow 0,
$$
where $\frac{1}{p-\eta}-\frac{1}{q-\varepsilon}=\alpha$.
\end{proof}
\vskip+0.2cm

{\bf Corollary 3.1.} {\em Let $\theta >0$ and let $1<p<\infty$. Suppose that $0<\alpha <1/p$. We set $q=\frac{p}{1-\alpha p}$. Then $I_{\alpha}$ is bounded from $L^{p), \theta_1}([0,1])$ to $L^{q), \theta_2}([0,1])$ provided that $\theta_2 > (1+\alpha q)\theta_1$.}
\vskip+0.2cm

{\em Proof} follows immediately from Theorem 3.1 (in the unweighted case $w(x)\equiv \; const$) and $(1.1)$. $\;\;\;\;\;\;\; \Box $

\section{One-sided potentials}

\

In this section we show that the unboudedness result in grand
Lebesgue spaces is also true for the one--sided potentials:

$$
(R_{\alpha}f)(x)=\int\limits_{0}^{x}\frac{f(t)}{(x-t)^{1-\alpha}}dt,\;\;\;x\in
(0,1);
$$
and
$$
(W_{\alpha}f)(x)=\int\limits_{x}^{1}\frac{f(t)}{(t-x)^{1-\alpha}}dt,\;\;\;x\in
(0,1),
$$
where $0<\alpha<1$. In particular, we claim that $R_{\alpha}$ and $W_{\alpha}$ are not bounded
from $L^{p),\theta_{1}}$ to $L^{q),\theta_{2}}$, where
$q=\frac{p}{1-\alpha p}$, $1<p<\infty$, $\theta_{1}$, $\theta_{2}>0$, $\theta_{2}<\frac{\theta_{1}q}{p}$.
Indeed, let us show the result first for
$R_{\alpha}$.

Suppose the contrary:
$$
\|R_{\alpha}f\|_{L^{q),\theta_{2}}([0,1])}\leq c
\|f\|_{L^{p),\theta_{1}}([0,1])}, \;\;\;\;\;
\theta_{2}<\frac{\theta_{1}q}{p},  \eqno{(4.1)}
$$
where $c$ does not depend on $f$. Let
$f_{n}(x)=\chi_{(0,1/2n)}(x)$ in (4.1). Then taking the following
inequality
$$
(R_{\alpha}f_n)(x)\geq\int\limits_{0}^{\frac{1}{2n}}\frac{1}{(x-t)^{1-\alpha}}dt
\geq \big(\frac{1}{2n}\big)^{\alpha},\;\;\;\;
x\in\left(\frac{1}{2n},\frac{1}{n}\right), \eqno{(4.2)}
$$
into account, (4.1) yields that
$$
(2n)^{-\alpha}\bigg\|\chi_{\left(\frac{1}{2n},\frac{1}{n}\right)}\bigg\|_{L^{q),\theta_{2}}([0,1])}
\leq c\bigg\|\chi_{(0,1/2n)}\bigg\|_{L^{p),\theta_{1}}([0,1])}.
\eqno{(4.3)}
$$
Now we choose $\varepsilon_{n}$ positive number so that
$$
\sup\limits_{0<\varepsilon\leq
p-1}\left(\varepsilon^{\theta_{1}}\frac{1}{2n}\right)^{\frac{1}{p-\varepsilon}}
=\left(\varepsilon_{n}^{\theta_{1}}\frac{1}{2n}
\right)^{\frac{1}{p-\varepsilon_{n}}}. \eqno{(4.4)}
$$
We now observe that $\lim\limits_{n\rightarrow
0}\varepsilon_{n}=0$ (see the proof of Theorem 2.1 for the similar
arguments). Choose now $\eta_{n}$ so that
$$
\alpha=\frac{1}{p}-\frac{1}{q}=\frac{1}{p-\varepsilon_{n}}-\frac{1}{q-\eta_{n}}.
$$
Hence,

$$ \eta_{n}=q-\frac{p-\varepsilon_{n}}{1-\alpha(p-\varepsilon_{n})}.
\eqno{(4.5)} $$
By (4.3)-(4.5) we conclude that

$$
(2n)^{-\alpha}\eta_{n}^{\frac{\theta_{2}}{q-\eta_{n}}}\left(\frac{1}{2n}\right)^{\frac{1}{q-\eta_{n}}}
\leq c
\varepsilon_{n}^{\frac{\theta_{1}}{p-\varepsilon_{n}}}(2n)^{-1/(p-\varepsilon_{n})}.
\eqno{(4.6)}
$$
From (4.6) we have that

$$
\eta_{n}^{\frac{\theta_{2}}{q-\eta_{n}}}\varepsilon_{n}^{-\frac{\theta_{1}}{p-\varepsilon_{n}}}
\leq c_{p},\;\;\text{for all}\;  n\in N \eqno{(4.7)}
$$
because

$$
\frac{1}{2}\leq\left(\frac{1}{2}\right)^{\frac{1}{p-\varepsilon_{n}}}\leq
\left(\frac{1}{2}\right)^{\frac{1}{p}},
$$
$$
\frac{1}{2}\leq\left(\frac{1}{2}\right)^{\frac{1}{q-\eta_{n}}}
\leq \left(\frac{1}{2}\right)^{\frac{1}{q}}.
$$

Now (4.5) yields
$$
\left[\frac{q-\frac{p-\varepsilon_{n}}{1-\alpha(p-\varepsilon_{n})}}{\varepsilon_{n}}\right]
^{\frac{\theta_{2}}{p-\varepsilon_{n}}-\alpha\theta_{2}}\cdot\varepsilon_{n}^{-\frac{\theta_{1}}
{p-\varepsilon_{n}}+\frac{\theta_{2}}{p-\varepsilon_{n}}-\alpha\theta_{2}}\leq
c_{p}.
$$
Hence,
$$
\left[\frac{q-\frac{p-\varepsilon_{n}}{1-\alpha(p-\varepsilon_{n})}}{\varepsilon_{n}}\right]
^{\frac{\theta_{2}}{p-\varepsilon_{n}}-\alpha\theta_{2}}\varepsilon_{n}^{\frac{\theta_{2}-\theta_{1}}
{p-\varepsilon_{n}}-\alpha\theta_{2}}\leq c_{p},
$$
which is impossible, because
$\lim\limits_{n\rightarrow\infty}\varepsilon_{n}^{\frac{\theta_{2}-\theta_{1}}
{p-\varepsilon_{n}}-\alpha\theta_{2}}=\infty$ (recall that
$\frac{\theta_{2}-\theta_{1}} {p}-\alpha\theta_{2}=
\frac{\theta_2}{q}- \frac{\theta_1}{p}<0$).

Analogously, we have  that $W_{\alpha}$ is not bounded from
$L^{p),\theta_{1}}$ to $L^{q),\theta_{2}}$. This follows from the
inequalities
$$
(W_{\alpha})(x)\geq\int\limits_{x}^{1-\frac{1}{3n}}\frac{f(t)}{(t-x)^{1-\alpha}}dt\geq
\left(\frac{2}{3n}\right)^{\alpha-1}\cdot\frac{1}{6n}=c_{\alpha}n^{-\alpha},\;\;\;x\in
\left(1-\frac{1}{n},1-\frac{1}{2n}\right),
$$
where $f(t)= \chi_{(1-\frac{1}{2n}, 1-\frac{1}{3n})}(t)$.  Hence,
$$
c_{\alpha}n^{-\alpha}\bigg\|\chi_{\left(1-\frac{1}{n},1-\frac{1}{2n}\right)}\bigg\|_{L^{q),\theta_{2}}([0,1])}
\leq
c\bigg\|\chi_{\left(1-\frac{1}{2n},1-\frac{1}{3n}\right)}\bigg\|_{L^{p),\theta_{1}}([0,1])}.
$$
Choosing now $\varepsilon_{n}$ so that
$$
\left[\varepsilon_{n}^{\theta_{1}}\frac{1}{6n}\right]^{\frac{1}{p-\varepsilon_{n}}}=
\sup\limits_{0<\varepsilon_{n}\leq
p-1}\left[\varepsilon_{n}^{\theta_{1}}\frac{1}{6n}\right]^{\frac{1}{p-\varepsilon}},
\;\;\; 0<\varepsilon_{n}\leq p-1, $$ and observing that
$\lim\limits_{n\rightarrow\infty}\varepsilon_{n}=0$ (see the proof
of Theorem 2.1 for the similar arguments) we find that the
conclusion similar to the case of $R_{\alpha}$ is valid.

\subsection{Conclusions and Remarks}

\setcounter{equation}{0}

\

Let $0<\alpha<1$ and let $I_{\alpha}$, $R_{\alpha}$, $W_{\alpha}$
be potential operators defined above. In the sequel  we denote by
$T_{\alpha}$ one of these operators.


\vskip+0.2cm

{\bf Corollary 5.1.} {\em   Let $1<p<\infty$ and let
$0<\alpha<1/p$. We set $q=\frac{p}{1-\alpha p}$. Suppose that
$\theta_{1}$ and $\theta_{2}$ be positive numbers. Then:

$\rm{(i)}$\;If $\theta_{2}<(1+\alpha q)\theta_{1}$, then
$T_{\alpha}$ is not bounded from $L^{p),\theta_{1}}$ to
$L^{q),\theta_{2}}$.

$\rm{(ii)}$\;If $\theta_{2}\geq(1+\alpha q)\theta_{1}$, then
$T_{\alpha}$ is bounded from $L^{p),\theta_{1}}$ to
$L^{q),\theta_{2}}$.}

\vskip+0.2cm

{\em Remark} 5.1. There is a function $f$ from $L^{p)}\backslash
L^{p}$ such that $T_{\alpha}f\in L^{q)}\backslash L^{q}$.

Indeed, let $f(t)=t^{-\frac{1}{p}}$, $t\in (0,1)$. Then $f\in
L^{p)}\backslash L^{p}$. On the other hand, (see e. g.
\cite{SaKiMa}), $T_{\alpha}f\approx t^{-\frac{1}{q}}$. Hence
$T_{\alpha}f\in L^{q)}\backslash L^{q}$.



\section*{Acknowledgements}

The author expresses his gratitude to Professor V. Kokilashvili
for careful reading of the manuscript and valuable comments and
suggestions. The author thanks also Professor L. Ephremidze for
the discussions regarding the proof of Theorem 2.1.


\vskip+1cm

\vskip+0.5cm

 Author's Address:

\

A. Meskhi: \

A. Razmadze Mathematical Institute, M. Aleksidze St.,  Tbilisi
0193, Georgia

Second Address: Department of Mathematics,  Faculty of Informatics
and Control Systems, Georgian Technical University, 77, Kostava
St., Tbilisi, Georgia.

e-mail:  meskhi@rmi.acnet.ge \vskip+0.5cm

\end{document}